\documentclass[12pt,twoside]{article}
\usepackage[english]{babel}
\usepackage[latin1]{inputenc}
\usepackage{amsmath}
\usepackage{amssymb,amsfonts}
\usepackage{graphicx}
\usepackage{times,amssymb,amscd}
\newtheorem{thm}{Theorem}[section]

\newtheorem{lem}[thm]{Lemma}
\newtheorem{prop}[thm]{Proposition}
\newtheorem{defn}[thm]{Definition}
\newtheorem{rem}[thm]{Remark}
\numberwithin{equation}{section}

\newcommand{\bA}{\mathbf{A}}

\newcommand{\bE}{\mathbf{E}}

\newcommand{\bH}{\mathbf{H}}

\newcommand{\bL}{\mathbf{L}}

\newcommand{\bR}{\mathbf{R}}
\newcommand{\bS}{\mathbf{S}}
\newcommand{\bV}{\mathbf{V}}

\newcommand{\be}{\mathbf{e}}

\newcommand{\br}{\mathbf{r}}
\newcommand{\bx}{\mathbf{x}}

\newcommand{\bT}{\mathbf{T}}

\newcommand{\BV}{\boldsymbol{V}}

\newcommand{\Be}{\boldsymbol{e}}
\newcommand{\Bu}{\boldsymbol{u}}
\newcommand{\Bv}{\boldsymbol{v}}

\newcommand{\cP}{\mathcal{P}}
\newcommand{\cS}{\mathcal{S}}

\newcommand{\EUC}{\mathbf E^3}
\newcommand{\SPH}{\bS^3}
\newcommand{\HYP}{\bH^3}
\newcommand{\SXR}{\bS^2\!\times\!\bR}
\newcommand{\HXR}{\bH^2\!\times\!\bR}
\newcommand{\SLR}{\widetilde{\bS\bL_2\bR}}
\newcommand{\NIL}{\mathbf{Nil}}
\newcommand{\SOL}{\mathbf{Sol}}

\begin{document}
\pagestyle{myheadings}
\markboth{\centerline{Angéla Vr\'anics and Jen\H o Szirmai}}
{Lattice coverings by congruent translation balls $\dots$}
\title
{Lattice coverings by congruent translation balls using translation-like bisector surfaces in $\NIL$ geometry
}

\author{Angéla Vr\'anics and Jen\H o Szirmai \\
\normalsize Budapest University of Technology and \\
\normalsize Economics, Institute of Mathematics, \\
\normalsize Department of Geometry \\
\normalsize Budapest, P. O. Box: 91, H-1521 \\
\normalsize szirmai@math.bme.hu
\date{\normalsize{\today}}}

\maketitle
\begin{abstract}
In this paper we study the $\NIL$ geometry that is one of the eight homogeneous Thurston 3-geomet\-ri\-es.

We determine the equation of the translation-like bisector surface
of any two points. We prove, that the isosceles property of a translation triangle is not equivalent to two angles of the triangle being equal and that
the triangle inequalities do not remain valid for translation triangles in general.
We develop a method to determine the centre and the radius of the circumscribed translation sphere of a given {\it translation tetrahedron}.

A further aim of this paper is to study lattice-like coverings with congruent translation balls in $\NIL$
space. We introduce the notion of the density of the considered coverings
and give upper estimate to it using the radius amd the volume of the circumscribed translation sphere of a given {\it translation tetrahedron}. 
The found minimal upper bound density of the translation ball coverings $\Delta \approx 1.42783$.
In our work we will use for computations and visualizations the projective model of $\NIL$ described by E. Moln\'ar in \cite{M97}.
\end{abstract}

\section{Introduction} \label{section1}
The basic problems in the classical theory of packings and coverings, the development of which was strongly influenced by the geometry of numbers 
and by crystallography, are the determination of the densest packing and the thinnest covering with congruent copies of a given body. At present the body is a ball and now
we consider the lattice-like covering problem with congruent translation balls in $\NIL$ space.

These questions related to the theory of the Dirichlet-Voronoi cells (brifly $D-V$ cells).
In $3$-dimensional spaces of constant curvature the $D-V$ cells are widely investigated, but in the further Thurston geometries $\SXR$, $\HXR$,
$\NIL$, $\SOL$, $\SLR$ there are few results in this topic. Let $X$ be one of the above five geometries and $\Gamma$ is one of its discrete isometry groups.
Moreover, we distinguish two distance function types: $d^g$ is the usual geodesic distance function and $d^t$ is the translation distance function (see Section 3).
Therefore, we obtain two types of the $D-V$ cells regarding the two distance functions.

The firs step to get the $D-V$ cell of a given point set of $X$ is the determination of the translation or geodesic-like bisector (or equidistant) surface
of two arbitrary points of $X$ because these surface types contain the faces of $D-V$ cells. 

In \cite{PSchSz10}, \cite{PSchSz11}, \cite{PSchSz11-1} we studied the geodesic-like equidistant surfaces in $\SXR$, $\HXR$ and $\NIL$ geometries, and in \cite{Sz17-1} we discussed the 
translation-like bisector surfaces in $\SOL$ geometry, but there are no results concerning the translation-like equidistant surfaces in $\NIL$ and $\SLR$ geometries. 

In the Thurston spaces can be introduced in a natural way (see \cite{M97}) translations mapping each point to any point.
Consider a unit vector at the origin. Translations, postulated at the
beginning carry this vector to any point by its tangent mapping. If a curve $t\rightarrow (x(t),y(t),z(t))$ has just the translated
vector as tangent vector in each point, then the  curve is
called a {\it translation curve}. This assumption leads to a system of first order differential equations, thus translation
curves are simpler than geodesics and differ from them in $\NIL$, $\SLR$ and $\SOL$ geometries. In $\EUC$, $\SPH$, $\HYP$, $\SXR$ and $\HXR$ geometries the translation and geodesic
curves coincide with each other.

Therefore, the translation curves also play an important role in $\NIL$, $\SLR$ and $\SOL$ geometries and often seem to be more natural in these geometries,
than their geodesic lines.  

{\it In this paper} we study the translation-like bisector surface of any two points in $\NIL$ geometry, determine its equation and visualize them.
The translation-like bisector surfaces play an important role in the construction of the $D-V$ cells because their faces lie on bisector surfaces. The $D-V$-cells
are relevant in the study of tilings, ball packing and ball covering. E.g. if the point set is the orbit of a point - generated by
a discrete isometry group of $\NIL$ - then we obtain a monohedral $D-V$ cell decomposition (tiling) of the considered space and it is interesting to examine its
optimal ball packing and covering (see \cite{Sz13-1}, \cite{Sz13-2}).

Moreover, we prove, that the isosceles property of a translation triangle is not equivalent to two angles of the triangle being equal and that
the triangle inequalities do not remain valid for translation triangles in general.

Using the above bisector surfaces we develop a procedure to determine the centre and the radius of the circumscribed translation sphere of an arbitrary $\NIL$ tetrahedron.
This is useful to determine the least dense ball covering radius of a given periodic polyhedral $\NIL$ tiling because the tiling can be decomposed into tetrahedra.
Applying the above procedure we determine the minimal covering density of some lattice types and thus we give an upper bound of the lattice-like covering 
density related to the most important lattice parameter $k=1$.
\section{On $\NIL$ geometry}
$\NIL$ geometry can be derived from the famous real matrix group $\mathbf{L(\mathbb{R})}$ discovered by Werner {Heisenberg}. The left (row-column) 
multiplication of Heisenberg matrices
     \begin{equation}
     \begin{gathered}
     \begin{pmatrix}
         1&x&z \\
         0&1&y \\
         0&0&1 \\
       \end{pmatrix}
       \begin{pmatrix}
         1&a&c \\
         0&1&b \\
         0&0&1 \\
       \end{pmatrix}
       =\begin{pmatrix}
         1&a+x&c+xb+z \\
         0&1&b+y \\
         0&0&1 \\
       \end{pmatrix}
      \end{gathered} \tag{2.1}
     \end{equation}
defines "translations" $\mathbf{L}(\mathbb{R})= \{(x,y,z): x,~y,~z\in \mathbb{R} \}$ 
on the points of $\NIL= \{(a,b,c):a,~b,~c \in \mathbb{R}\}$. 
These translations are not commutative in general. The matrices $\mathbf{K}(z) \vartriangleleft \mathbf{L}$ of the form
     \begin{equation}
     \begin{gathered}
       \mathbf{K}(z) \ni
       \begin{pmatrix}
         1&0&z \\
         0&1&0 \\
         0&0&1 \\
       \end{pmatrix}
       \mapsto (0,0,z)  
      \end{gathered}\tag{2.2}
     \end{equation} 
constitute the one parametric centre, i.e. each of its elements commutes with all elements of $\mathbf{L}$. 
The elements of $\mathbf{K}$ are called {\it fibre translations}. $\NIL$ geometry of the Heisenberg group can be projectively 
(affinely) interpreted by "right translations" 
on points as the matrix formula 
     \begin{equation}
     \begin{gathered}
       (1;a,b,c) \to (1;a,b,c)
       \begin{pmatrix}
         1&x&y&z \\
         0&1&0&0 \\
         0&0&1&x \\
         0&0&0&1 \\
       \end{pmatrix}
       =(1;x+a,y+b,z+bx+c) 
      \end{gathered} \tag{2.3}
     \end{equation} 
shows, according to (2.1). Here we consider $\mathbf{L}$ as projective collineation group with right actions in homogeneous coordinates.
We will use the Cartesian homogeneous coordinate simplex $E_0(\be_0)$,$E_1^{\infty}(\be_1)$,$E_2^{\infty}(\be_2)$,
$E_3^{\infty}(\be_3), \ (\{\be_i\}\subset \bV^4$ \ $\text{with the unit point}$ $E(\be = \be_0 + \be_1 + \be_2 + \be_3 ))$ 
which is distinguished by an origin $E_0$ and by the ideal points of coordinate axes, respectively. 
Moreover, $\mathbf{y}=c\bx$ with $0<c\in \mathbb{R}$ (or $c\in\mathbb{R}\setminus\{0\})$
defines a point $(\bx)=(\mathbf{y})$ of the projective 3-sphere $\mathcal{P} \mathcal{S}^3$ (or that of the projective space $\cP^3$ where opposite rays
$(\bx)$ and $(-\bx)$ are identified). 
The dual system $\{(\Be^i)\}, \ (\{\Be^i\}\subset \BV_4)$, with $\be_i\Be^j=\delta_i^j$ (the Kronecker symbol), describes the simplex planes, especially the plane at infinity 
$(\Be^0)=E_1^{\infty}E_2^{\infty}E_3^{\infty}$, and generally, $\Bv=\Bu\frac{1}{c}$ defines a plane $(\Bu)=(\Bv)$ of $\cP \cS^3$
(or that of $\cP^3$). Thus $0=\bx\Bu=\mathbf{y}\Bv$ defines the incidence of point $(\bx)=(\mathbf{y})$ and plane
$(\Bu)=(\Bv)$, as $(\bx) \text{I} (\Bu)$ also denotes it. Thus {\bf Nil} can be visualized in the affine 3-space $\bA^3$
(so in $\bE^3$) as well \cite{MSz06}.

In this context E. Moln\'ar \cite{M97} has derived the well-known infinitesimal arc-length square invariant under translations $\bL$ at any point of $\NIL$ as follows
\begin{equation}
   \begin{gathered}
      (dx)^2+(dy)^2+(-xdy+dz)^2=\\
      (dx)^2+(1+x^2)(dy)^2-2x(dy)(dz)+(dz)^2=:(ds)^2
       \end{gathered} \tag{2.4}
     \end{equation}
The translation group $\mathbf{L}$ defined by formula (2.3) can be extended to a larger group $\mathbf{G}$ of collineations,
preserving the fibres, that will be equivalent to the (orientation preserving) isometry group of $\NIL$. 

In \cite{M06} E.~Moln\'ar has shown that 
a rotation through angle $\omega$
about the $z$-axis at the origin, as isometry of $\NIL$, keeping invariant the Riemann
metric everywhere, will be a quadratic mapping in $x,y$ to $z$-image $\overline{z}$ as follows:
     \begin{equation}
     \begin{gathered}
       \mathcal{M}=\br(O,\omega):(1;x,y,z) \to (1;\overline{x},\overline{y},\overline{z}); \\ 
       \overline{x}=x\cos{\omega}-y\sin{\omega}, \ \ \overline{y}=x\sin{\omega}+y\cos{\omega}, \\
       \overline{z}=z-\frac{1}{2}xy+\frac{1}{4}(x^2-y^2)\sin{2\omega}+\frac{1}{2}xy\cos{2\omega}.
      \end{gathered} \tag{2.5}
     \end{equation}
This rotation formula $\mathcal{M}$, however, is conjugate by the quadratic mapping $\alpha$ to the linear rotation $\Omega$ in (1.7) as follows
     \begin{equation}
     \begin{gathered}
       \alpha^{-1}: \ \ (1;x,y,z) \stackrel{\alpha^{-1}}{\longrightarrow} (1; x',y',z')=(1;x,y,z-\frac{1}{2}xy) \ \ \text{to} \\
       \Omega: \ \ (1;x',y',z') \stackrel{\Omega}{\longrightarrow} (1;x",y",z")=(1;x',y',z')
       \begin{pmatrix}
         1&0&0&0 \\
         0&\cos{\omega}&\sin{\omega}&0 \\
         0&-\sin{\omega}&\cos{\omega}&0 \\
         0&0&0&1 \\
       \end{pmatrix}, \\
       \text{with} \ \ \alpha: (1;x",y",z") \stackrel{\alpha}{\longrightarrow}  (1; \overline{x}, \overline{y},\overline{z})=(1; x",y",z"+\frac{1}{2}x"y").
      \end{gathered} \tag{2.6}
     \end{equation}
This quadratic conjugacy modifies the $\NIL$ translations in (2.3), as well. Now a translation with $(X,Y,Z)$ in (2.3) instead of $(x,y,z)$ will be changed 
by the above conjugacy to the translation 
     \begin{equation}
     \begin{gathered}
        (1;x,y,z) \longrightarrow (1; \overline{x}, \overline{y},\overline{z})=(1; x,y,z)
        \begin{pmatrix}
         1&X&Y&Z-\frac{1}{2}XY \\
         0&1&0&-\frac{1}{2}Y \\
         0&0&1&\frac{1}{2}X \\
         0&0&0&1 \\
       \end{pmatrix}, \\
             \end{gathered} \tag{2.7}
     \end{equation}
that is again an affine collineation.
\subsection{Translation curves and balls}
We consider a $\NIL$ curve $(1,x(t), y(t), z(t) )$ with a given starting tangent vector at the origin $O=E_0=(1,0,0,0)$
\begin{equation}
   \begin{gathered}
      u=\dot{x}(0),\ v=\dot{y}(0), \ w=\dot{z}(0).
       \end{gathered} \tag{2.8}
     \end{equation}
For a translation curve let its tangent vector at the point $(1,x(t), y(t), z(t) )$ be defined by the matrix (2.3) 
with the following equation:
\begin{equation}
     \begin{gathered}
     (0,u,v,w)
     \begin{pmatrix}
         1&x(t)&y(t)&z(t) \\
         0&1&0&0 \\
         0&0&1&x(t) \\
         0&0&0&1 \\
       \end{pmatrix}
       =(0,\dot{x}(t),\dot{y}(t),\dot{z}(t)).
       \end{gathered} \tag{2.9}
     \end{equation}
Thus, the {\it translation curves} in $\NIL$ geometry (see  \cite{MoSzi10}, \cite{MSz06} \cite{MSzV}) are defined by the above first order differential equation system 
$\dot{x}(t)=u, \ \dot{y}(t)=v,  \ \dot{z}(t)=v \cdot x(t)+w,$ whose solution is the following: 
\begin{equation}
   \begin{gathered}
       x(t)=u t, \ y(t)=v t,  \ z(t)=\frac{1}{2}uvt^2+wt.
       \end{gathered} \tag{2.10}
\end{equation}
We assume that the starting point of a translation curve is the origin, because we can transform a curve into an 
arbitrary starting point by translation (2.3), moreover, unit initial velocity translation can be assumed by "geographic" parameters
$\phi$ and $\theta$:
\begin{equation}
\begin{gathered}
        x(0)=y(0)=z(0)=0; \\ \ u=\dot{x}(0)=\cos{\theta} \cos{\phi}, \ \ v=\dot{y}(0)=\cos{\theta} \sin{\phi}, \ \ w=\dot{z}(0)=\sin{\theta}; \\ 
        - \pi \leq \phi \leq \pi, \ -\frac{\pi}{2} \leq \theta \leq \frac{\pi}{2}. \tag{2.11}
\end{gathered}
\end{equation}
\begin{defn}
The translation distance $d^t(P_1,P_2)$ between the points $P_1$ and $P_2$ is defined by the arc length of the above translation curve 
from $P_1$ to $P_2$.
\end{defn}
\begin{defn} The sphere of radius $r >0$ with centre at the origin, (denoted by $S^t_O(r)$), with the usual longitude and altitude parameters 
$\phi$ and $\theta$, respectively by (2.11), is specified by the following equations:
\begin{equation}
\begin{gathered}
        S_O^t(r): \left\{ \begin{array}{ll} 
        x(\phi,\theta)=r \cos{\theta} \cos{\phi}, \\
        y(\phi,\theta)=r \cos{\theta} \sin{\phi}, \\
        z(\phi,\theta)=\frac{r^2}{2} \cos^2{\theta} \cos{\phi} \sin{\phi}+ r \sin{\theta}.
        \end{array} \right.
        \tag{2.12}
\end{gathered}
\end{equation}
\end{defn}
\begin{defn}
 The body of the translation sphere of centre $O$ and of radius $r$ in the $\NIL$ space is called translation ball, denoted by $B^t_{O}(r)$,
 i.e. $Q \in B^t_{O}(r)$ iff $0 \leq d^t(O,Q) \leq r$.
\end{defn}
\begin{rem}
The translation sphere is a simply connected 
 surface without selfintersection in $\NIL$ space for any radius $0<r\in \mathbb{R}$.
\end{rem}

We obtained in \cite{Sz12-1} the volume formula of the translation ball $B^t_{O}(r)$ of radius $r$ by (2.4), (2.5) and (2.12):
\begin{thm}
The volume of a translation ball of radius $r$ is the same as that of an Euclidean one: 
\begin{equation}
Vol(B^t_{O}(r))=\frac{4}{3}r^3 \pi. \tag{2.13}
\end{equation}
\end{thm}
The convexity of the translation ball play an important role in the discussion of the ball covering therefore we recall the following Theorem from the paper \cite{Sz12-1}.
\begin{thm}
A translation $\NIL$ ball $B^t(S^t(r))$ is convex in the affine-Euclidean sense in our model if and only if $r \in [0, 2]$. 
\end{thm}
\subsection{The discrete translation group $L(\mathbb{Z}$,~{\it{k}})}
We consider the $\NIL$ translations defined in (2.1) and (2.3) and choose first two non-commuting translations 
\begin{equation}
\begin{gathered}
\tau_1=
\begin{pmatrix}
1&t_1^1&t_1^2&t_1^3 \\
0&1&0&0 \\
0&0&1&t_1^1 \\
0&0&0&1 \\
\end{pmatrix} \ \text{and} \ 
\tau_2=
\begin{pmatrix}
1&t_2^1&t_2^2&t_2^3 \\
0&1&0&0 \\
0&0&1&t_2^1 \\
0&0&0&1 \\
\end{pmatrix},
\end{gathered} \tag{2.14}
\end{equation}
now with upper indices for the coordinate variables.
Second, we define the translation $(\tau_3)^k$,  $(k $ $\in $ $\mathbb{N}\setminus\{0\} \ \text{$k$ is fixed natural exponent)},$ by the following commutator:
\begin{equation}
(\tau_3)^k=\tau_2^{-1}\tau_1^{-1}\tau_2\tau_1=
\begin{pmatrix}
1 & 0 & 0 & -t_2^1 t_1^2+t_1^1 t_2^2 \\
0&1&0&0\\
0&0&1&0\\
0&0&0&1
\end{pmatrix}, \ \ \text{and so $\tau_3$} \tag{2.15}
\end{equation}
is also defined. If we take integers as coefficients for $\tau_1,\tau_2,\tau_3,$ then we generate the discrete group 
$\langle \tau_1,\tau_2, \tau_3 \rangle$,  denoted by $\mathbf{L}(\tau_1,\tau_2,k )$ or by $\mathbf{L}(\mathbb{Z},k)$. Here
$\mathbb{Z}$ refers to the integers.

We know (see e.g. \cite{S} and \cite{Sz07}) that the orbit space $\NIL /\mathbf{L}(\mathbb{Z},k)$ is a compact manifold, i.e. a $\NIL$ space form.
\begin{defn}
The $\NIL$ point lattice $\Gamma_P (\tau_1,\tau_2,k )$ is a discrete orbit of point $P$ in the $\NIL$ space under group $\mathbf{L}(\tau_1,\tau_2,k )$= 
$\mathbf{L}(\mathbb{Z},k)$ 
with an arbitrary starting point $P$ for every fixed $k \in \mathbb{N}\setminus \{0\}.$
\end{defn}
\begin{rem}
For simplicity we have chosen the origin as starting point, by the homogeneity of $\NIL$. 
\end{rem}
\begin{rem}
We may assume in the following that $t_1^2=0$, i.e. the image of the origin by the translation $\tau_1$ lies on the plane $[x,z]$.
\end{rem}

We consider by (2.14-15) a fundamental "parallelepiped complex" (see \cite{Sz12-1})
$$\widetilde{\mathcal{F}(k)}=OT_1T_2T_3T_{12}T_{21}T_{23}T_{213}T_{13},\ \  \text{(see Fig.~1 for $k=1,2$)}$$
in the Euclidean sense, which is determined by translations $\tau_1,\tau_2,\tau_3$.
The images of $\widetilde{\mathcal{F}(k)}$ under $\mathbf{L}(\mathbb{Z},k)$ fill $\NIL$ without gap. Overlaps occur only on the boundary.
\begin{figure}[ht]
\centering
\includegraphics[width=6cm]{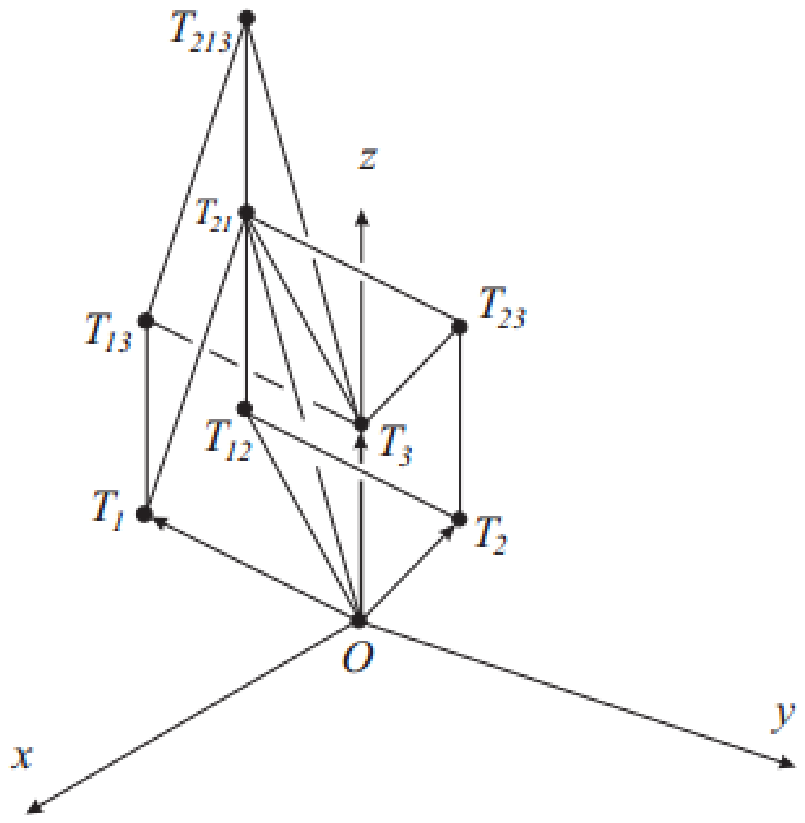} ~ \includegraphics[width=4.5cm]{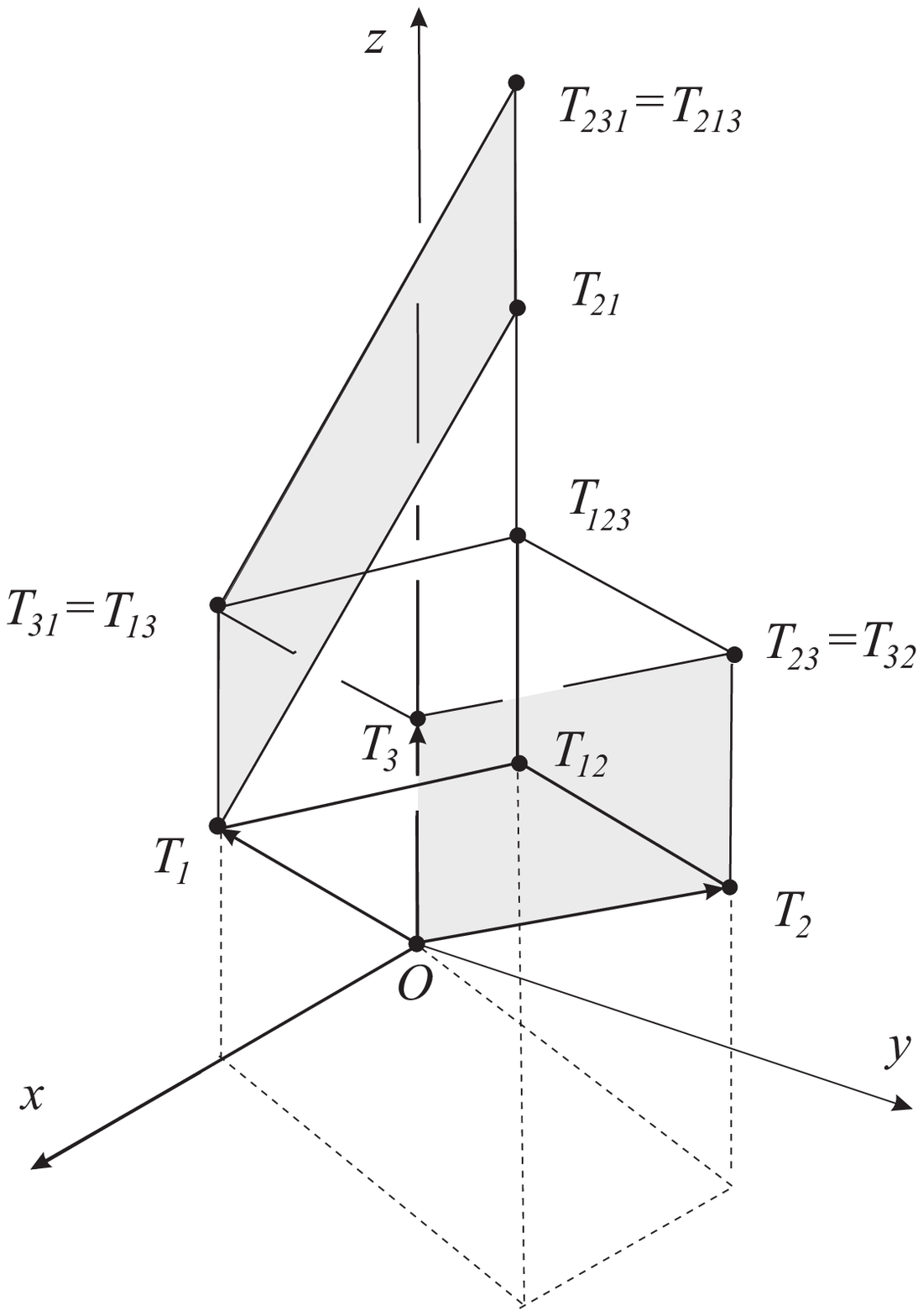}
\caption{The $\NIL$ parallelepipeds $\widetilde{\mathcal{F}(1)}$ (left) and $\widetilde{\mathcal{F}(2)}$ (right)}
\label{}
\end{figure}
{\it Analogously to the Euclidean integer lattice and parallelepiped, $\widetilde{\mathcal{F}(k)} \ \ (k \in \mathbb{N}\setminus \{0 \})$ can be called a 
$\NIL$ parallelepiped, endowed by face pairing, as the upper $\sim$ hints to it.}

$\widetilde{\mathcal{F}(k)}$ is a {\it fundamental domain} of $\mathbf{L}(\mathbb{Z},k)$. We need only its interior for its volume.
The homogeneous coordinates of the vertices of $\widetilde{\mathcal{F}(k)}$ can be determined in our affine model by the translations
(2.14-15) with the parameters $t_i^j, \ i\in \{1,2\}, \ j \in \{1,2,3\}$ (see (2.16) and Fig.~2). 
\begin{equation}
\begin{gathered}
T_1(1,t_1^1,0,t_1^3), \ T_2(1,t_2^1,t_2^2,t_2^3), \ T_3(1,0,0,\frac{t_1^1 t_2^2}{k}), \\ 
T_{13}(1,t_1^1,0,\frac{t_1^1 t_2^2}{k}+t_1^3), \ T_{12}(1,t_1^1+t_2^1,t_2^2,t_2^3+t_1^3), \\ 
T_{21}(1,t_1^1+t_2^1,t_2^2,t_1^1 t_2^2+t_1^3+t_2^3), \ 
T_{23}(1,t_2^1,t_2^2,t_2^3+\frac{t_1^1 t_2^2}{k}), \\ T_{213}=T_{231}(1,t_1^1+t_2^1,t_2^2, (k+1) \frac{t_1^1 t_2^2}{k}+t_1^3+t_2^3). \tag{2.16}
\end{gathered}
\end{equation}
In \cite{Sz07} we have determined the volume of the $\NIL$ parallelepiped $\widetilde{\mathcal{F}(1)}$. Analogously to that we get the volume
formula of $\widetilde{\mathcal{F}(k)}$ ($k \in \mathbb{N}$) by the usual method:   
\begin{equation}
\int\!\!\int\!\!\int_{\widetilde{\mathcal{F}(k)}} {\sqrt{det(g_{ij})}} ~ \mathrm{d} x \mathrm{d} y \mathrm{d} z = 
Vol(\widetilde{\mathcal{F}(k)})=\frac{1}{k}\int_0^{t_2^2} \int_0^{t_1^1} ~|t_1^1 \cdot t_2^2| ~ \mathrm{d}{x}\mathrm{d}{y}=\frac{(t_1^1 \cdot t_2^2)^2}{k}. \tag{2.17}
\end{equation}
If the parameter $k$ is given, from this formula it can be seen that the volume of a $\NIL$ parallelepiped depends on two parameters, 
i.e. on its projection into the $[x,y]$ plane.
\section{Translation-like bisector surfaces}
Our further goals are to examine and visualize the Dirichlet-Voronoi cells and the packing and covering problems 
of $\NIL$ geometry. In order to study the above questions have to determine 
the "faces" of the $D-V$ cells that are
parts of bisector (or equidistant) surfaces of given point pairs. 
The definition below comes naturally:
\begin{defn}
The equidistant surface $\cS_{P_1P_2}$ of two arbitrary points $P_1,P_2 \in \NIL$ consists of all points $P'\in \NIL$,
for which $d^t(P_1,P')=d^t(P',P_2)$.
\end{defn}
It can be assumed by the homogeneity of $\NIL$ that the starting point of a 
given translation curve segment is $E_0=P_1=(1,0,0,0)$ and 
the other endpoint will be given by its homogeneous coordinates $P_2=(1,a,b,c)$. 
We consider the translation curve segment $t_{P_1P_2}$ and determine its
parameters $(\phi,\theta,r)$ expressed by the real coordinates $a$, $b$, $c$ of $P_2$. 
We obtain directly by equation system (2.12) the following:
\begin{lem}
\begin{enumerate}
\item Let $(1,a,b,c)$ $(a,b \in \mathbb{R} \setminus \{0\},~c\in \mathbb{R})$ be the homogeneous 
coordinates of the point $P \in \NIL$. The parameters of the
corresponding translation curve $t_{E_0P}$ are the following
\begin{equation}
\begin{gathered}
\phi=\mathrm{arccot}\Big(\frac{a}{b}\Big),~ 
\theta=\mathrm{arccot}\Big(\frac{\sqrt{a^{2}+b^{2}}}{c-\frac{ab}{2}}\Big),~
r=\Big|\frac{c-\frac{ab}{2}}{\sin{\theta}}\Big|.
\end{gathered} \tag{3.1}
\end{equation}
\item Let $(1,a,0,c)$ $(a,c \in \mathbb{R} \setminus \{0\})$ be the homogeneous 
coordinates of the point $P \in \NIL$. The parameters of the
corresponding translation curve $t_{E_0P}$ are the following
\begin{equation}
\begin{gathered}
\phi=\pi \cdot n ,~ (n\in\{0,1\}),~ 
\theta= \mathrm{arccot}\Big(\frac{a}{c}\Big),~
r=\Big|\frac{a}{\cos{\theta}}\Big|.
\end{gathered} \tag{3.2}
\end{equation}
\item Let $(1,a,0,0)$ $(a \in \mathbb{R}\setminus \{0\})$ 
be the homogeneous coordinates of the point $P \in \NIL$. The parameters of the 
corresponding translation curve $t_{E_0P}$ are the following
\begin{equation}
\begin{gathered}
\phi=\pi \cdot n ,~ (n\in\{0,1\}),~
\theta=\pi \cdot n ,~ (n\in\{0,1\}),~
r=|a|.
\end{gathered} \tag{3.3}
\end{equation}
\item Let $(1,0,b,0)$ $(b \in \mathbb{R}\setminus \{0\})$ 
be the homogeneous coordinates of the point $P \in \NIL$. The parameters of the 
corresponding translation curve $t_{E_0P}$ are the following
\begin{equation}
\begin{gathered}
\phi=\pm \frac{\pi}{2},~
\theta=\pi \cdot n ,~ (n\in\{0,1\}),~
r=|b|.
\end{gathered} \tag{3.4}
\end{equation}
\item Let $(1,0,0,c)$ $(c \in \mathbb{R}\setminus \{0\})$ 
be the homogeneous coordinates of the point $P \in \NIL$. The parameters of the 
corresponding translation curve $t_{E_0P}$ are the following
\begin{equation}
\begin{gathered}
\theta=\pm \frac{\pi}{2},~
r=|c|.~ ~ \square
\end{gathered} \tag{3.5}
\end{equation}
\end{enumerate}
\end{lem}
{\it In order to determine the translation-like bisector surface 
$\cS_{P_1P_2}(x,y,z)$ of two given point $E_0=P_1=(1,0,0,0)$ and $P_2=(1,a,b,c)$
we define the {translation} $\bT_{P_2}$ as elements of the isometry group of 
$\NIL$, that
maps the origin $E_0$ onto $P_2$} (see Fig.~2), 
moreover let $P_3=(1,x,y,z)$ a point in $\NIL$ space.

This isometrie $\bT_{P_2}$ and its inverse (up to a positive determinant factor) can be given by:
\begin{equation}
\bT_{P_2}=
\begin{pmatrix}
1 & a & b & c \\
0 & 1 & 0 & 0 \\
0 & 0 & 1 & a \\
0 & 0 & 0 & 1
\end{pmatrix} , ~ ~ ~
\bT_{P_2}^{-1}=
\begin{pmatrix}
1 & -a & -b & ab-c \\
0 & 1 & 0 & 0 \\
0 & 0 & 1 & -a \\
0 & 0 & 0 & 1
\end{pmatrix} , \tag{3.6}
\end{equation}
and the images $\bT^{-1}_{P_2}(P_i)$ of points $P_i$ $(i \in \{1,2,3\})$ are the following (see also Fig.~2):
\begin{equation}
\begin{gathered}
\bT^{-1}_{P_2}(P_1=E_0)=P_1^2=(1,-a,-b,ab-c),~ \bT^{-1}_{P_2}(P_2)=E_0=(1,0,0,0), \\
\bT^{-1}_{P_2}(P_3)=P_3^2=(1,(x-a),(y-b),a(b-y)-c). \tag{3.7}
\end{gathered}
\end{equation}
It is clear that $P_3=(1,x,y,z) \in \cS_{P_1P_2} ~ \text{iff} ~ d^t(P_1,P_3)=d^t(P_3,P_2) \Rightarrow d^t(P_1,P_3)=d^t(E_0,P_3^2)$ where
$P_3^2=\bT^{-1}_{P_2}(P_3)$ (see (3.6), (3.7)).

This method leads to
\begin{lem}
The equation of the equidistant surface $\cS_{P_1P_2}(x,y,z)$ of 
two points $P_1=(1,0,0,0)$ and $P_2=(1,a,b,c)$ in $\NIL$ space (see Fig.~2,3):
\begin{enumerate}
\item $a,b,c \ne 0,$
\begin{equation}\label{bis1}
\begin{gathered}
z=\frac{1}{4} \Big( \frac{8 x \left(a^2+b^2\right)-4 \left(a^3-ab+4 bc\right)}{a (b (a+x)-
ay-2c)}-\\- \frac{b (a (a+x)+8)}{a}+y (a+2 x)+2 c \Big),
\end{gathered} \tag{3.8}
\end{equation}
\item $a,b \ne 0,~c=0$
\begin{equation}\label{bis1}
\begin{gathered}
z= -\frac{a^2 \left(b^2-2 b y+y^2+4\right)+2 a x \left(b^2-2 b y+y^2-4\right)}{4
   (a (b-y)+b x)}-\\-\frac{b\left(x^2+4\right) (b-2 y)}{4(a (b-y)+b x)},
\end{gathered} \tag{3.9}
\end{equation}
\item $a,c \ne 0,~b=0$
\begin{equation}\label{bis1}
\begin{gathered}
z= \frac{a^2 \left(y^2+4\right)+2 a \left(2 c y+x \left(y^2-4\right)\right)+
4 c (c+x y)}{4 a y+8 c} ,
\end{gathered} \tag{3.10}
\end{equation}
\item $b,c \ne 0,~a=0$
\begin{equation}\label{bis1}
\begin{gathered}
z= \frac{b^2 \left(x^2+4\right)-2 b \left(2 c x+\left(x^2+4\right) y\right)+
4 c (c+x y)}{8 c-4 b x},
\end{gathered} \tag{3.11}
\end{equation}
\item $b,c=0,~a\ne0$
\begin{equation}\label{bis1}
\begin{gathered}
z=\frac{a \left(y^2+4\right)+2 x \left(y^2-4\right)}{4 y}).
\end{gathered} \tag{3.12}
\end{equation}
\item $a,b=0,~b\ne0$
\begin{equation}\label{bis1}
\begin{gathered}
z=-\frac{\left(x^2+4\right) (b-2 y)}{4 x}.
\end{gathered} \tag{3.13}
\end{equation}
\item $a,b=0,~c\ne0$
\begin{equation}\label{bis1}
\begin{gathered}
z=\frac{1}{2} (c+x y).~\square
\end{gathered} \tag{3.14}
\end{equation}
\end{enumerate}
\end{lem}
\begin{figure}[ht]
\centering
\includegraphics[width=13cm]{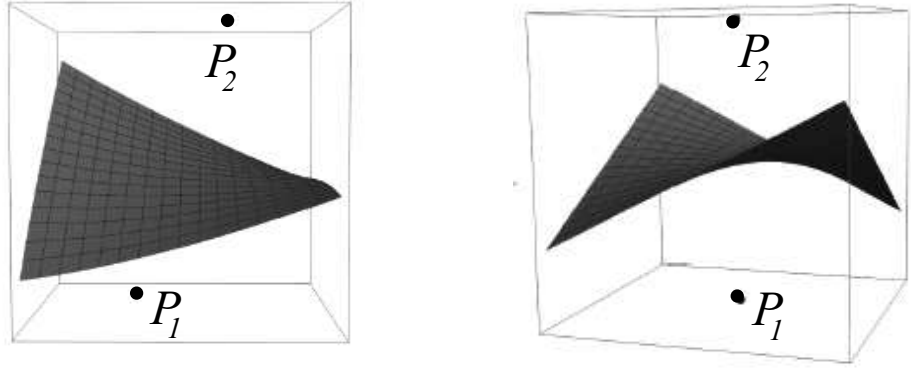} 
\caption{Translation-like bisectors (equidistant surfaces) of point pairs 
$(P_1,P_2)$ with coordinates $((1,0,0,0), (1,1/2,1/2,3/2))$ (left) 
and $((1,0,0,0), (1,0,0,2))$ (right) }
\label{pic:surf1}
\end{figure}
\subsection{On isosceles and equilateral translation triangles}
We consider $3$ points $A_1$, $A_2$, $A_3$ in the projective model of $\NIL$ space.
The {\it translation segments} connecting the points $A_i$ and $A_j$
$(i<j,~i,j,k \in \{1,2,3\})$ are called sides of the {\it translation triangle}
$A_1A_2A_3$. The length of its side $a_k$
$(k\in \{1,2,3\})$ is the translation distance $d^t(A_i,A_j)$ between the vertices 
$A_i$ and $A_j$ $(i<j,~i,j,k \in \{1,2,3\}, k \ne i,j$).

Similarly to the Euclidean geometry we can define the notions of isosceles and 
equilateral translation triangles.
\begin{figure}[ht]
\centering
\includegraphics[width=.95\linewidth]{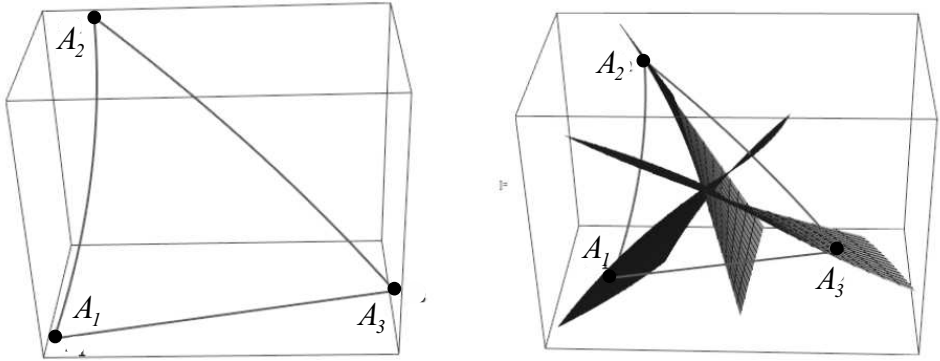} 
\caption{Equidistant surfaces of the "edges" of the equilateral triangle 
$A_1A_2A_3$ where the coordinates of the vertices $A_1(1,0,0,0)$, $A_2(1,0.8,0.5,-0.131662)$, 
$A_3(1,0.2,-0.058102,-0.983882)$).}
\label{pic:surf1}
\end{figure}
An isosceles translation triangle is a triangle with (at least) two equal sides and a triangle
with all sides equal is called an equilateral translation triangle (see Fig.~3) in the $\NIL$ space.

We note here, that if in a translation triangle $A_1A_2A_3$ e.g. $a_1=a_2$ then 
the bisector surface $\cS_{A_1A_2}$ contains the vertex $A_3$ (see Fig.~3).

In the Euclidean space the isosceles property of a triangle is equivalent to 
two angles of the triangle being equal
therefore has both two equal sides and two equal angles. 
An equilateral triangle is a special case of an isosceles triangle 
having not just two, but all three sides and angles equal.
\begin{prop}
The isosceles property of a translation triangle is not equivalent to two angles of the triangle being equal in the $\NIL$ space.
\end{prop}
{\bf Proof:}~
The missing coordinates $y^3$ and $z^3$ of the vertices $A_1=E_0=(1,0,0,0,0)$, $A_2=(1,x^2=1,y^2=1/2,z^2=-3/4)$ and
$A_3=(1,x^3=0,y^3,z^3)$ can be determined by the equation system $d^t(A_1,A_2)=d^t(A_1,A_3)=d^t(A_2,A_3)$. We get the following coordinates:
$y^3 \approx -0.6164636$, $z^3 \approx -1.367469$ where $(a_3=d^t(A_1,A_2)=a_2=d^t(A_1,A_3)=a_1=d^t(A_2,A_3) = 1.5)$.

The {\it interior angles} of translation triangles are denoted at the vertex $A_i$ by $\omega_i$ $(i\in\{1,2,3\})$.
We note here that the angle of two intersecting translation curves depends on the orientation of their tangent vectors. 

{\it In order to determine the interior angles of a translation triangle $A_1A_2A_3$
and its interior angle sum $\sum_{i=1}^3(\omega_i)$,}
we apply the method (we do not discuss here) developed in \cite{Sz17} using the infinitesimal arc-lenght square of $\NIL$ geometry (see (2.4)).

Our method (see \cite{Sz17}) provide the following results:
$$
\omega_1\approx 1.08063\text{,  }\omega_2\approx 0.84167\text{,  }\omega_3\approx 1.22186,~\sum_{i=1}^3(\omega_i) \approx 3.14416 > \pi \approx 3.14159.
$$
From the above results follows the statement. We note here, that if the vertices of the translation triangle lie in the $[x,y]$ plane than the
Euclidean isosceles property true in the $\NIL$ geometry, as well. ~$\square$

Using the above methods we obtain the following
\begin{lem}
The triangle inequalities do not remain valid for translation triangles in general.
\end{lem}
{\bf Proof:}~
We consider the translation triangle $A_1A_2A_3$ where $A_1=(1,0,0,0)$, 
$A_2=(1,-1,3,1)$, 
$A_3=(1,1/4,1/2,1/2)$. 
We obtain directly by equation systems (3.1-5) (see Lemma 3.2 and \cite{Sz17}) 
the lengths of the translation segments $A_iA_j$ $(i,j \in \{1,2,3\}$, $i<j)$:
\begin{equation}
\begin{gathered}
d^t(A_1,A_2) \approx 4.03113,~d^t(A_1,A_3) \approx 0.70986, ~ d^t(A_2,A_3) \approx 3.14307, \\ \text{therefore}
~ d^t(A_2,A_3)+d^t(A_1,A_3) < d^t(A_1,A_2). ~\square \notag
\end{gathered}
\end{equation}
\subsection{The locus of all points equidistant from three given points}
A point is said to be equidistant from a set of objects if the distances between that point and each object in the set are equal.
Here we study that case where the objects are vertices of a $\NIL$ translation triangle $A_1A_2A_3$ and determine the locus of all points 
that are equidistant from $A_1$, $A_2$ and $A_3$.

We consider $3$ points $A_1$, $A_2$, $A_3$ that do not all lie in the same translation curve in the projective model of $\NIL$ space.
The {\it translation segments} connecting the points $A_i$ and $A_j$ $(i<j,~i,j,k \in \{1,2,3\}, k \ne i,j$) are called sides of the {\it translation triangle} $A_1A_2A_3$.
The locus of all points that are equidistant from the vertices $A_1$, $A_2$ and $A_3$ is denoted by $\mathcal{C}$.

In the previous section we determined the equation of translation-like bisector (equidistant) surface to any two points in the $\NIL$ space.
It is clear, that all points on the locus $\mathcal{C}$ must lie on the equidistant surfaces $\cS_{A_iA_j}$, $(i<j,~i,j \in \{1,2,3\})$ therefore
$\mathcal{C}=\mathcal{S}_{A_1A_2} \cap \mathcal{S}_{A_1A_3}$ and the coordinates of each of the points of that locus and only those points must satisfy the corresponding
equations of Lemma 3.3. Thus, the non-empty point set $\mathcal{C}$ can be determined and can be visualized for any given translation triangle (see Fig.~4 and 5).
In the Fig.~4 we describe the translation triangle $A_1A_2A_3$ with vertices $A_1=(1,0,0,0)$, $A_2=(1,0,0,1)$, $A_3=(1,0,1,0)$ with the equidistant surfaces 
$$
\mathcal{S}_{A_1A_2}:~z=\frac{1}{8} (4 x y+4), ~ ~ \mathcal{S}_{A_2A_3}:~z=\frac{2 x y^2-8 x+y^2+4 y+4}{4 y}
$$
of edges $A_1A_2$ and $A_2A_3$ and their intersection $\mathcal{C}=\mathcal{S}_{A_1A_2} \cap \mathcal{S}_{A_2A_3}$.
\begin{figure}[ht]
\centering
\includegraphics[width=12cm]{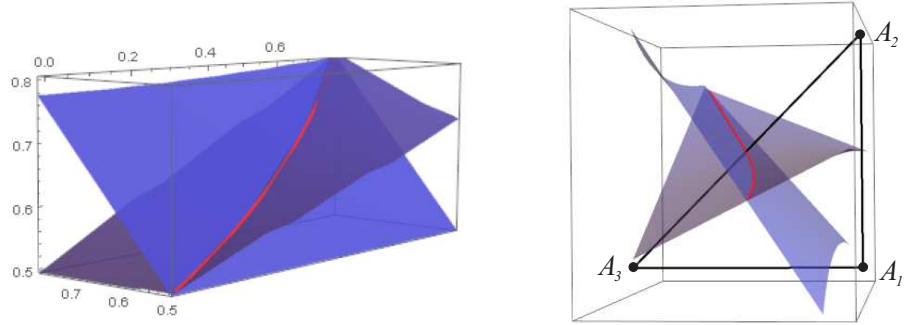}
\caption{Translation triangle with vertices $A_1=(1,0,0,0)$, $A_2=(1,0,0,1)$, $A_3=(1,0,1,0)$ with translation-like bisector surfaces $\cS_{A_1A_2}$ and $\cS_{A_2A_3}$ 
and a part of the locus $\mathcal{C}=\cS_{A_1A_2} \cap \cS_{A_1A_3}$ of all points equidistant from three given points $A_1$, $A_2$, $A_3$.}
\label{}
\end{figure}
If the vertices of the translation triangle $A_1A_2A_3$ lie in e.g. coordinate plane $[y,z]$ or $[x,z]$ we obtain the following lemmas:
\begin{lem}
If the vertices of a translation triangle $A_1A_2A_3$ lie on the $[y,z]$ plane $A_1=(1,0,0,0)$, $A_2=(1,0,b_2,b_3)$, $A_3=(1,0,c_2,c_3)$ 
($b_2\neq 0$, $b_3\neq 0$, $c_2\neq 0$, $c_3\neq 0$) then
the parametric equation $(x\in \mathbb{R})$ of $\mathcal{C}$ is the following (see Lemma 3.3 and Fig.~5):
\begin{equation*}
\mathcal{C}(x):~\left\{x,\frac{f}{16 (b_2 c_3-b_3 c_2)},\frac{g}{32 (b_3
  c_2-b_2 c_3)}\right\}
\end{equation*} {\text where}
\begin{equation*}
\begin{gathered}
f=-2 b_3 \left(-2 c_2 x (b_2 x+2 c_3)+4 c_3 (b_2 x+c_3)+c_2^2 \left(x^2+4\right)\right)+\\+b_2 \left(b_2
   \left(x^2+4\right)(2c_3-c_2 x)  
   +x \left(c_2^2 \left(x^2+4\right)-4 c_2 c_3 x+4 c_3^2\right)\right)+\\+b_3^2 (8 c_3-4 c_2
   x),
\end{gathered}
\end{equation*} ~ {\text and}~
\begin{equation*}
\begin{gathered}g=b_2^2 (x^2+4)  (c_2 (x^2+4)-2 c_3 x)-b_2 (4 c_2 x
   (x^2+4) (b_3-c_3)+\\+4 c_3 (c_3 (x^2+4)
   -2 b_3 x^2)+c_2^2 (x^2+4)^2)+2 b_3 (2
   b_3 (c_2(x^2+4)-2 c_3 x)+ \\ +x (c_2^2 (x^2+4)-4 c_2 c_3 x+4 c_3^2)).
   \end{gathered}
   \end{equation*}
\end{lem}

\begin{lem}
If the vertices of a translation triangle $A_1A_2A_3$ lie on the $[x,z]$ plane $A_1=(1,0,0,0)$, $A_2=(1,b_1,0,b_3)$, $A_3=(1,c_1,0,c_3)$ 
($b_1\neq 0$, $b_3\neq 0$, $c_1\neq 0$, $c_3\neq 0$) then
the parametric equation $(y\in \mathbb{R})$ of $\mathcal{C}$ is the following (see Lemma 3.3 and Fig.5):
\begin{equation*}
\mathcal{C}(y):~\left\{\frac{f}{16 (b_1 c_3-b_3 c_1)},y,\frac{g}{16 (b_3 c_1-b_1 c_3)}\right\}
\end{equation*} where 
\begin{equation*}
\begin{gathered}
f=-2 b_3 \left(-2 c_1 y (b_1 y-2 c_3)+4 c_3 (c_3-b_1 y)+c_1^2 \left(y^2+4\right)\right)+\\
+b_1 \left(b_1
   \left(y^2+4\right) (c_1 y+2 c_3)-y \left(c_1^2 \left(y^2+4\right)+4 c_1 c_3 y+4 c_3^2\right)\right)+\\+4 b_2^2 (c_1 y+2
   c_3),
\end{gathered}
\end{equation*} and 
\begin{equation*}
\begin{gathered}
g=-b_1^2 c_1 y^3-4 b_1^2 c_1 y-2 b_1^2 c_3 y^2-8 b_1^2 c_3-4
   b_1 b_3 c_1 y^2-8 b_1 b_3 c_3 y+\\+b_1 c_1^2 y^3 +4 b_1 c_1^2 y+4 b_1 c_1 c_3 y^2 +4 b_1
   c_3^2 y-4 b_3^2 c_1 y-8 b_3^2 c_3+\\+2 b_3 c_1^2 y^2+8 b_3 c_1^2+8 b_3 c_1 c_3 y+8 b_3
   c_3^2.
\end{gathered}
\end{equation*}
\end{lem}
\begin{figure}[ht]
\centering
\includegraphics[width=14cm]{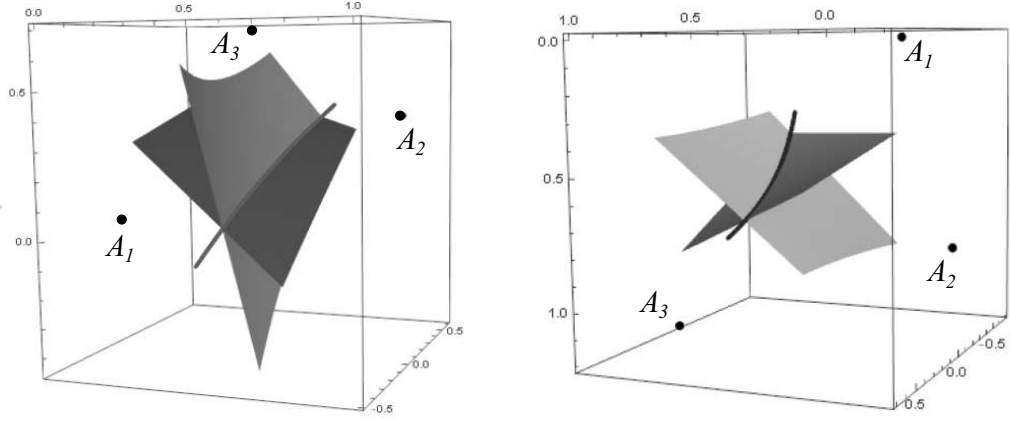}
\caption{Translation triangle with vertices $A_1=E_0=(1,0,0,0)$, $A_2=(1,1/2,0,7/10)$, $A_3=(1,1,0,2/5)$ with translation-like bisector surfaces $\cS_{A_1A_2}$ and $\cS_{A_1A_3}$ (left)
and Translation triangle with vertices $A_1=E_0=(1,0,0,0)$, $A_2=(1,0,-1/5,4/5)$, $A_3=(1,0,1,6/5)$ with translation-like bisector surfaces $\cS_{A_1A_2}$ and $\cS_{A_1A_3}$ (right).}
\label{}
\end{figure}
\subsection{Translation tetrahedra and their circumscribed spheres}
We consider $4$ points $A_1$, $A_2$, $A_3$, $A_4$ in the projective model of $\NIL$ space (see Section 2).
These points are the vertices of a {\it translation tetrahedron} in the $\NIL$ space if any two {\it translation segments} connecting the points $A_i$ and $A_j$
$(i<j,~i,j \in \{1,2,3,4\}$) do not have common inner points and any three vertices do not lie in a same translation curve.
Now, the translation segments $A_iA_j$ are called edges of the translation tetrahedron $A_1A_2A_3A_4$.
\begin{figure}[ht]
\centering
\includegraphics[width=12cm]{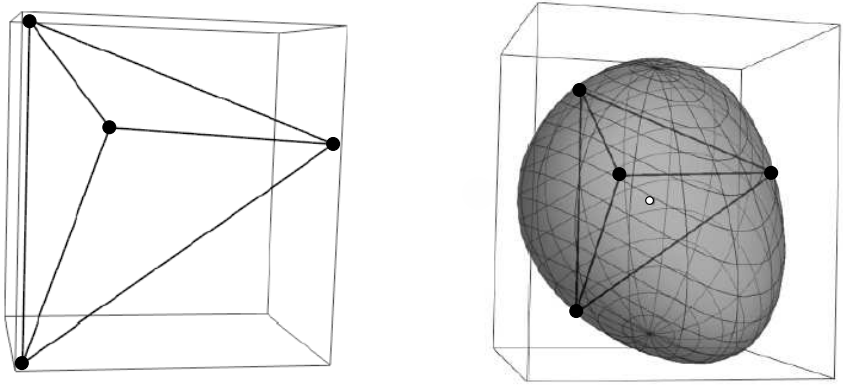}
\caption{Translation tetrahedron with vertices $A_1=(1,0,0,0)$, $A_2=(1,1.4,0,1)$, $A_3=(1,0.5,1,1)$, $A_4=(1,0,0,1.5)$
and its circumscibed sphere of radius $r \approx 0.92804$.}
\label{}
\end{figure}
The circumscribed sphere of a translation tetrahedron is a translation sphere (see Definition 2.2, (2.12)) that touches each of the tetrahedron's vertices.
As in the Euclidean case the radius
of a translation sphere circumscribed around a tetrahedron $T$ is called the circumradius of $T$, and the center point of this sphere is called the circumcenter of $T$.

\begin{lem}
For any translation tetrahedron there exists uniquely a translation sphere (called the circumsphere) on which all four vertices lie.
\end{lem}
\begin{figure}[ht]
\centering
\includegraphics[width=13cm]{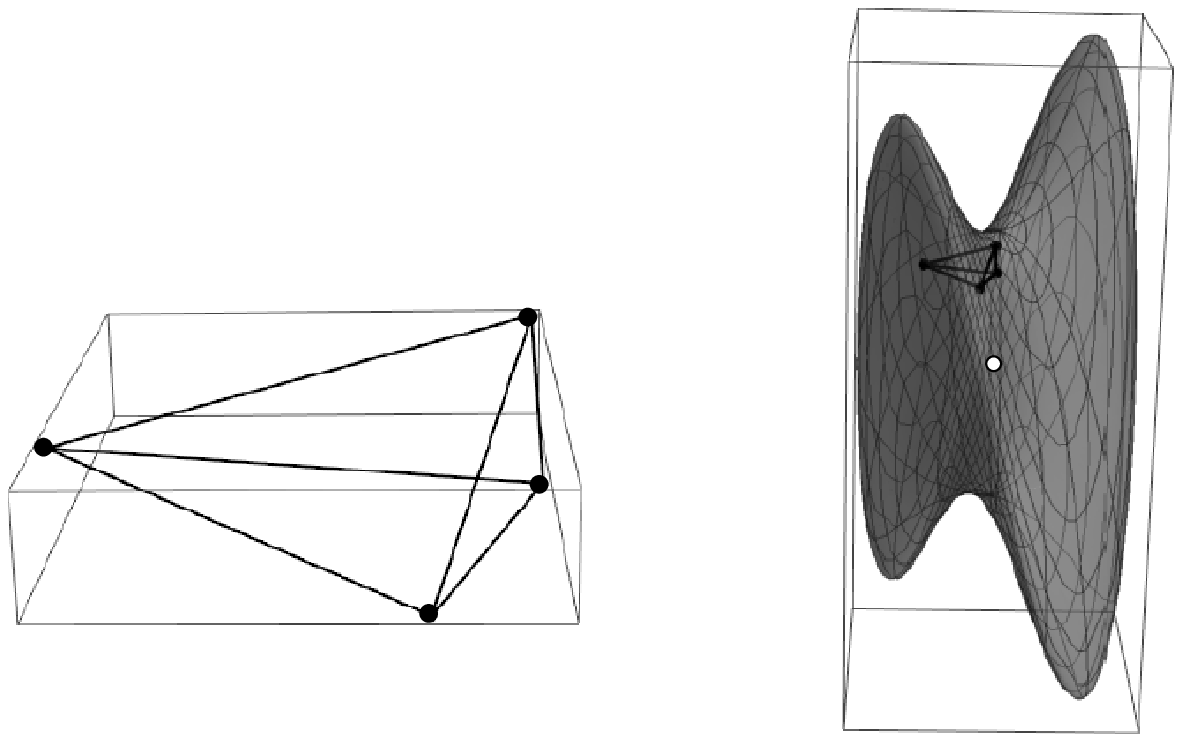}
\caption{Translation tetrahedron with vertices $A_1=(1,0,0,0)$, $A_2=(1,4,2,1)$, $A_3=(1,1,3,0)$, $A_4=(1,0,-2,1)$
and its circumscibed sphere of radius $r \approx 7.96825$.}
\label{}
\end{figure}
{\bf Proof:}~ The Lemma follows directly from the properties of the translation distance function (see Definition 2.1 and (2.12)).
The procedure to determine the radius and the circumcenter of a given translation tetrahedron is the folowing:

The circumcenter $C=(1,x,y,z)$ of a given translation tetrahedron $A_1A_2A_3A_4$ $(A_i=(1,x^i,y^i,z^i), ~ i \in \{1,2,3,4\})$
have to hold the following system of equation:
\begin{equation}
d^t(A_1,C)=d^t(A_2,C)=d^t(A_3,C)=d^t(A_4,C), \tag{3.15}
\end{equation}
therefore it lies on the translation-like bisector surfaces $\cS_{A_i,A_j}$ $(i<j,~i,j \in \{1,2,3,4\}$) which equations are determined in Lemma 3.3.
The coordinates $x,y,z$ of the circumcenter of the circumscribed sphere around the tetrahedron $A_1A_2A_3A_4$ are obtained by the system of equation
derived from the facts:
\begin{equation}
C \in \cS_{A_1A_2}, \cS_{A_1A_3}, \cS_{A_1A_4}. \tag{3.16}
\end{equation}
Finally, we get the circumradius $r$ as the translation distance e.g. $r=d^t(A_1,C)$.

We apply the above procedure to two tetrahedra determined their centres and the
radii of their circumscribed balls that are described in  Fig.~6 and 7.~ $\square$
\section{The lattice-like translation ball coverings}
In \cite{Sz13-1} we investigated the lattice-like {\it geodesic ball coverings} with congruent geodesic balls and in this section we study the similar problem
of the translation ball coverings.

In the following, we shall consider lattice coverings, each of them consisting of congruent translation balls of $\NIL$.  
Let $\mathcal{B}^c_\Gamma(R)$ denote a translation ball covering of $\NIL$ space with balls $B^c(R)$ of radius $R$ where their 
centres give rise to a $\NIL$ point lattice 
$\Gamma(\tau_1,\tau_2,k)$ ($k\in \mathbb{N}^+$). $\widetilde{\mathcal{F}(k)}$ 
is an arbitrary $\NIL$ parallelepiped of this lattice (see Section 2.2). 
The images of $\widetilde{\mathcal{F}(k)}$ by our discrete translation group $\bL(\tau_1,\tau_2,k)=\bL(\mathbb{Z},k)$
cover the $\NIL$ space without gap. 
\begin{rem}
In the $\NIL$ geometry, similarly to the Euclidean space $\bE^d, \ (d\ge 1)$, an arbitrary lattice $\Gamma$ gives a lattice-like covering of equal
balls if the radius $R$ of the balls is large enough. For the geodesic ball packings it is not true because the geodesic balls should have a radius $R \in [0,2\pi]$ (see \cite{Sz07}).
\end{rem}
If we start with a translation-like lattice covering $\mathcal{B}^c_\Gamma(R)$ and shrink the balls until they finally do not cover the space any more, then the 
minimal radius defines the least dense covering to a given lattice $\Gamma(\tau_1,\tau_2,k)$.
The thresfold value $R^c_\Gamma$ is called {\it the minimal covering radius} of the point lattice $\Gamma(\tau_1,\tau_2,k)$:
\begin{equation}
R^c_\Gamma:=\min \{R: ~ \text{where} ~  \mathcal{B}^c_\Gamma(R) \ \text{lattice covering by} ~ \Gamma(\tau_1,\tau_2,k) \}. \tag{4.1}
\end{equation}

For the density of the packing it is sufficient to relate the volume of the minimal covering ball
to that of the solid $\widetilde{\mathcal{F}(k)}$. 

Analogously to the Euclidean case it can be defined the density $\Delta(\mathcal{B}_\Gamma^c(R))$ of the lattice-like geodesic
ball covering $\mathcal{B}_\Gamma^c(R)$:
\begin{defn}
\begin{equation}
\Delta(\mathcal{B}_\Gamma^c(R)):=\frac{Vol({B}_\Gamma^c(R))}{Vol(\widetilde{\mathcal{F}(k)})}, \tag{4.2}
\end{equation}
and its minimum $\Delta(\mathcal{B}_\Gamma^c(R_\Gamma^c))$ for radius $R_\Gamma^c$ in (4.1).
\end{defn}
{\it The main problem is that to which lattice $\Gamma(\tau_1,\tau_2,k)$ belongs the optimal minimal density where $k\in\mathbb{N}^+$ is a given parameter.} 
\begin{equation}
\Delta_{opt}(\mathcal{B}^c)=
\inf_{\Gamma} \Big\{ \Delta(\mathcal{B}_\Gamma^c(R_\Gamma^c)) \Big\}. \tag{4.3}
\end{equation}
and $\Gamma_{opt}^c$ denotes any optimal lattice, if it exists at all.
\begin{rem}
The covering radius is the radius of the circumsphere of the lattice's Dirichlet-Voronoi cell i.e. the largest distance
between the midpoint and the vertices of its Dirichlet-Voronoi cell, whose description deserves separate studies (see \cite{SchM15}).
\end{rem}
{\it In the following we study the most important case related to parameter $k=1$.} 
\subsection{Method to determination of densest lattice-like translation ball covering of a given lattice}
We develop an algorithm to determine the lattice-like thinnest ball covering of a given lattice $\Gamma(\tau_1,\tau_2,1)$. 

The lattice is generated by the translations $\tau_1$ and $\tau_2$ where their coordinates in the model are $t_i^{j} \ (i=1,2; \ j=1,2,3)$ (see (2.16)). 

The $\NIL$ parallelepiped 
$\widetilde{\mathcal{F}(1)}=E_0T_1^{opt}T_2T_3T_{12}T_{21}T_{23}T_{213}T_{13}$ 
is a {\it fundamental domain} of $\bL(\mathbb{Z},1)$.
The homogeneous coordinates of its vertices can be derived from the coordinates of $\tau_1$ and $\tau_2$ (see Fig.~1 and (2.3) with (2.16)).
We examine the {\it minimal covering radius} $R^c$ to the given lattice $\Gamma(\tau_1,\tau_2,1)$.
\begin{equation}
R^c_\Gamma:=\min \{R: ~ \text{where} ~  \mathcal{B}^c_\Gamma(R) \ \text{lattice covering by} ~ \Gamma(\tau_1,\tau_2,1) \}. \notag
\end{equation}
It is sufficient to investigate such ball arrangements $\mathcal{B}^c_{\Gamma}(R)$ where the balls cover $\widetilde{\mathcal{F}(1)}$ .

From (2.14-16) follows, that the fundamental parallelepiped $\widetilde{\mathcal{F}(1)}$ can be decomposed into Euclidean tetrahedra
$\{E_0, T_1, T_2, T_3 \}$,  $\{T_3, T_1, T_{23}, T_{13} \}$, $\{T_3, T_1,$ $ T_{23}, T_2 \}$, 
$\{T_{12}, T_1, T_{23},$ $T_2 \}$, $\{T_{12}, T_1, $ $T_{23}, $ $T_{13} \}$, $\{T_{12}, T_{21}, 
T_{23}, T_{13} \}$ which fill it just once. 
The radius $R_i~(i=1,2 \dots 6)$ of each circumscribed ball to the above point sets can be determined by the procedure described in the previous section. 
It is clear, that the lattice-like ball arrangement $\mathcal{B}^c_\Gamma(R^{c}_\Gamma)$ of radius  $R^{c}_\Gamma = \max\{R_i\}$ cover the 
fundamental parallelepiped $\widetilde{\mathcal{F}(1)}$ and thus the $\NIL$ space if the translation ball of radius $R^{c}_\Gamma$ is convex in Euclidean sense i.e. 
$R^{c}_\Gamma \in [0,2]$ (see Theorem 2.6).
\subsubsection{Upper boud for the covering density}
To have a comparison, first we consider our optimal lattice-like arrangement $\mathcal{B}_\Gamma^{p}(R_{p})$ for the {\it conjectured 
densest lattice-like translation ball packing} 
in the $\NIL$ space (see \cite{Sz17}). These balls will be blown up to a covering. This optimal lattice is given in \cite{Sz12-1} with parameters 
\begin{equation}
\begin{gathered}
t_1^{1,p}\approx 1.31225; ~ t_1^{3,p}=\frac{t_{3}^{3,p}}{2}; ~ t_2^{1,p}\approx 0.65613; ~ t_2^{2,p} \approx 1.13644; \\  t_2^{3,p}\approx 1.11847;~
r_{p}\approx 0.74565;~ t_{3}^{3,p}=2r_{p}. 
\end{gathered} \tag{4.4}
\end{equation}
This packing can be generated by the translations $\Gamma^p(\tau_1^{p},\tau_2^{p},1)$
where $\tau_1^{p}$ and $\tau_2^{p}$ are given by the above coordinates $t_i^{j, p} \ i=1,2; \ j=1,2,3$ (see (4.3)). Thus we obtain the 
neigbouring balls around an arbitrary ball of the packing $\mathcal{B}_\Gamma^{p}(R_{\Gamma^p}^c)$ by the lattice $\Gamma^{p}(\tau_1^{p},\tau_2^{p},1)$. 
We have ball "columns" in $z$-direction and in regular hexagonal projection onto the $[x,y]$-plane.
\begin{figure}[ht]
\centering
\includegraphics[width=4.5cm]{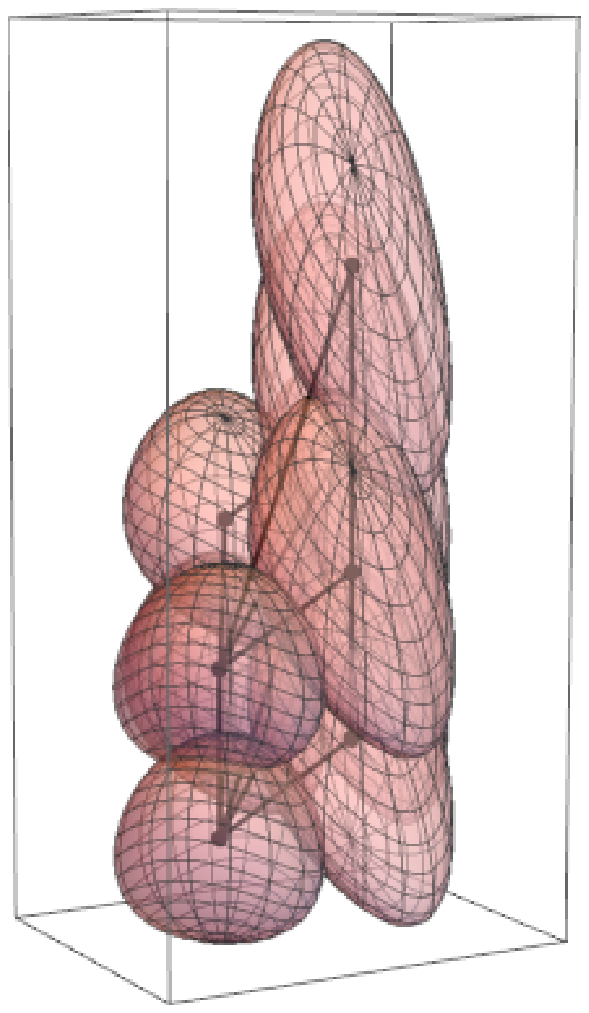}~ ~ ~ ~  \includegraphics[width=4cm]{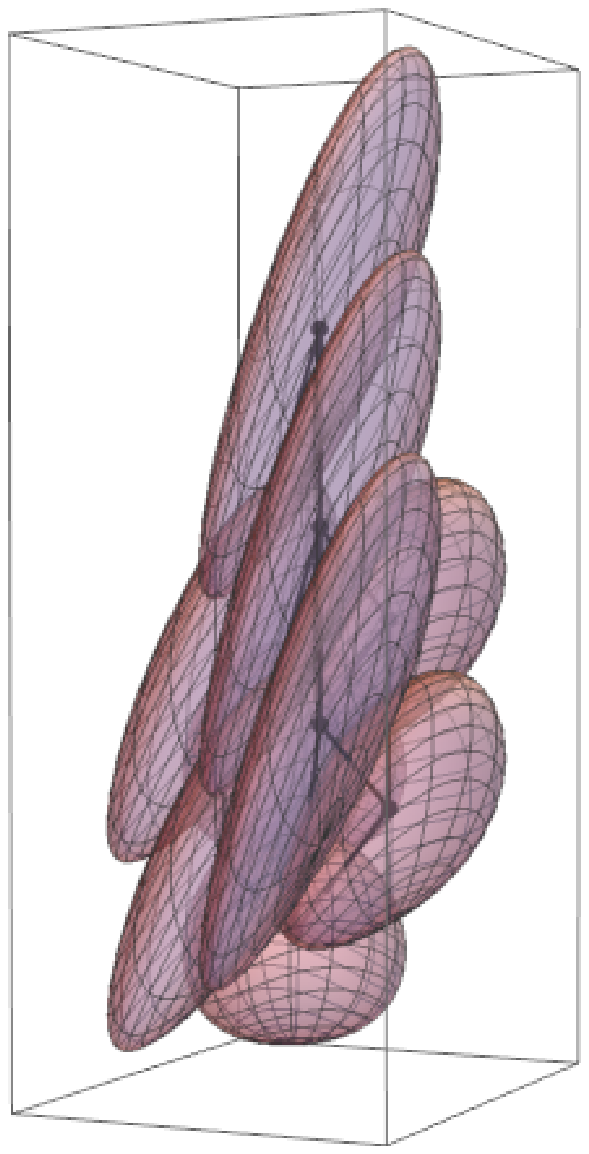}
\caption{Locally optimal lattice-like translation ball covering related to lattice $\Gamma^{p}(\tau_1^p,\tau_2^p,1)$.} 
\label{}
\end{figure}
From the structure of this lattice follows that in this case the corresponding lattice point sets $\{0, T_1^{p}, T_2^{p}, T_3^{p} \}$,  $\{T_3^{p}, $ $T_1^{p}, T_{23}^{p},$ 
$T_{13}^{p} \}$, $\{T_3^{p}, $ $T_1^{p},T_{23}^{p},$ $T_2^{p} \}$, 
$\{T_{12}^{p}, T_1^{p}, T_{23}^{p},$ $T_2^{p} \}$, $\{T_{12}^{p}, T_1^{p}, $ $T_{23}^{p}, $ $T_{13}^{p} \}$, $\{T_{12}^{p}, T_{21}^{p}, 
T_{23}^{p}, T_{13}^{p} \}$ are congruent by $\NIL$ isometries.
The radius $R_p$ of each circumscribed ball to the above point sets 
can be determined by the following system of equations:
\begin{equation}
\begin{gathered}
d^t(O,C)=d^t(C,T_3^{p})=d^t(C,T_1^{p})=d^t(C,T_2^{p}),
\end{gathered} \notag
\end{equation}
where $C(1,c^1,c^2,c^3)$ is the center of the circumscribed ball of the point set 
$\{E_0, T_1^{p}, $ $T_2^{p}, T_3^{p} \}$ ($d^t$ is the $\NIL$ translation distance, see Definition 2.1):
\begin{equation}
\begin{gathered}
c^1 \approx 0.45563, \ c^2 \approx  0.26306, \ c^3 \approx  0.80558, ~
R_{\Gamma^p}^c \approx 0.91257.
\end{gathered} \notag
\end{equation}
\begin{rem}
$C(1,c^1,c^2,c^3)$ is a vertex of the Dirichlet-Voronoi domain of the centre point $E_0$.
\end{rem}
$R_{\Gamma^p}^c \in [0,2]$ thus by Theorem 2.6 the ball of radius $R_{\Gamma^p}^c$ is convex in affin-Euclidean sense.  
Their circumscribed congruent $\NIL$ balls are convex thus they
cover the tetrahedra and so the ball arrangement $\mathcal{B}^c_{\Gamma^{p}}(R_{\Gamma^p}^c)$ cover the $\NIL$ space. 
Thus the radius $R_{\Gamma^p}^c$ of circumscribed ball 
give us the covering radius to the lattice $\Gamma^{p}$, indeed, and we get by (2.13), (2.17) and by the Definition
4.2 the following results:
\begin{equation}
\begin{gathered}
Vol(B(R_{\Gamma^p}^c)) \approx 3.18341, \ \ Vol(\widetilde{\mathcal{P}})=Vol(\widetilde{\mathcal{F}(1)}) \approx 2.22397, \\
\Delta(R_{\Gamma^p}^c,\tau^{p}_1,\tau^{p}_2,1):=\frac{Vol(\mathcal{B}_\Gamma(R_{\Gamma^p}^c))}{Vol(\widetilde{\mathcal{F}(1)})}
\approx 1.43141. \tag{4.5}
\end{gathered}
\end{equation}
\begin{rem}
The density of the least dense lattice-like ball covering in the the Euclidean space is 
$$
\Delta_{opt}(R^c_{opt},\tau^{c}_1,\tau^{c}_2,1) < \Delta_E=\frac{5 \sqrt{5} \pi}{24} \approx 1.46350.
$$
{This $\Delta_E$ attains for the so-called inner centred cubic lattice type of $\EUC$. That means a $\NIL$-lattice-ball-covering can be 
,,looser" than a Euclidean one.}
\end{rem}
Similarly to the above computations we can apply our method to any given $\NIL$ lattice. In the Table 1 we summarize the data of some locally optimal
lattice-like translation ball coverings:
\medbreak
\centerline{\vbox{
\halign{\strut\vrule\quad \hfil $#$ \hfil\quad\vrule
&\quad \hfil $#$ \hfil\quad\vrule &\quad \hfil $#$ \hfil\quad\vrule 
\cr
\noalign{\hrule}
\multispan3{\strut\vrule\hfill\bf Table 1 \hfill\vrule}%
\cr
\noalign{\hrule}
\noalign{\vskip2pt}
\noalign{\hrule}
\text{Lattice parameters} & R^c_\Gamma & \Delta^c_\Gamma   \cr
\noalign{\hrule}
t_i^j=1, ~ (i=1,2,~j=1,2,3) & \approx 0.88666 & \approx 2.91980  \cr
\noalign{\hrule}
\begin{gathered} t_1^1=t_1^{1,d},~ t_1^3=t_1^{3,d}, \\ t_2^1=t_2^{1,d},~ t_2^2=t_2^{2,d},~ t_2^3=t_2^{3,d} \end{gathered} & \approx 0.91257 & \approx 1.43141  \cr
\noalign{\hrule}
\begin{gathered}
t_1^1=1.3,~ t_1^3=0.74,\\ t_2^1=0.65,~ t_2^2=1.13,~ t_2^3=1.12
\end{gathered} & \approx 0.90406 & \approx 1.43429  \cr
\noalign{\hrule}
\begin{gathered}
t_1^1=1.29,~ t_1^3=0.74,\\ t_2^1=0.64,~ t_2^2=1.13,~ t_2^3=1.12
\end{gathered} & \approx 0.89997 & \approx 1.43692  \cr
\noalign{\hrule}
\begin{gathered} t_1^1=1.1,~ t_1^3=0.5,\\ t_2^1=0.5,~ t_2^2=1,~ t_2^3=1 \end{gathered} 
& \approx 0.77177 & \approx 1.59134  \cr
\noalign{\hrule}
\begin{gathered} t_1^1=1.1,~ t_1^3=0.5,\\ t_2^1=0.4,~ t_2^2=1,~ t_2^3=1 \end{gathered}
& \approx 0.78667 & \approx 1.68533  \cr
\noalign{\hrule}
\begin{gathered} t_1^1=1.31,~ t_1^3=0.74,\\ t_2^1=0.65,~ t_2^2=1.13,~ t_2^3=1.12 \end{gathered}
& \approx 0.90732 & \approx 1.42783  \cr
\noalign{\hrule}
}}}
\medbreak
From the above computations follows the following
\begin{thm}
The density of the least dense lattice-like translation ball covering is less or equal than the locally thinnest covering with congruent tranlation balls related to the lattice
$\Gamma^u(\tau_1^u,\tau_2^u,1)$ where the lattice is given by the parameters 
$t_1^1=1.31,~ t_1^3=0.74,~  t_2^1=0.65,~ t_2^2=1.13,~ t_2^3=1.12$ (see Fig.~9).
\begin{equation}
\Delta_{opt}(R^c_{opt},\tau^{c}_1,\tau^{c}_2,1) \le 
\Delta(R_{\Gamma^u}^c,\tau^{u}_1,\tau^{u}_2,1) \approx 1.42783 ~ 
(see ~ \text{Table} ~ 1 ~ and ~ Fig.~9). \notag
\end{equation}
\end{thm}
\begin{figure}[ht]
\centering
\includegraphics[width=4.5cm]{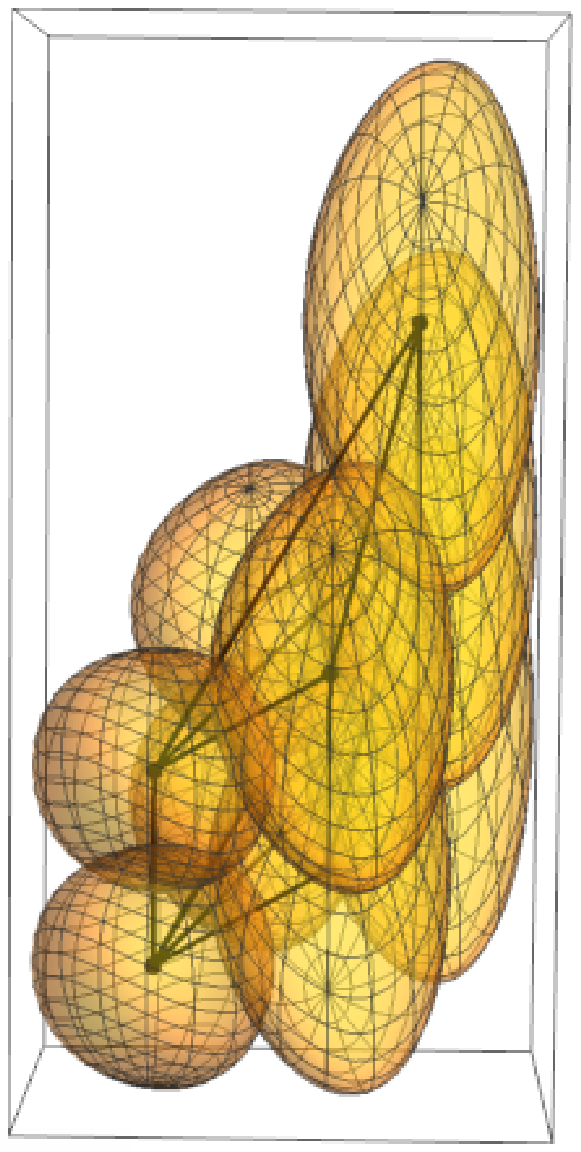} \quad \quad \includegraphics[width=5cm]{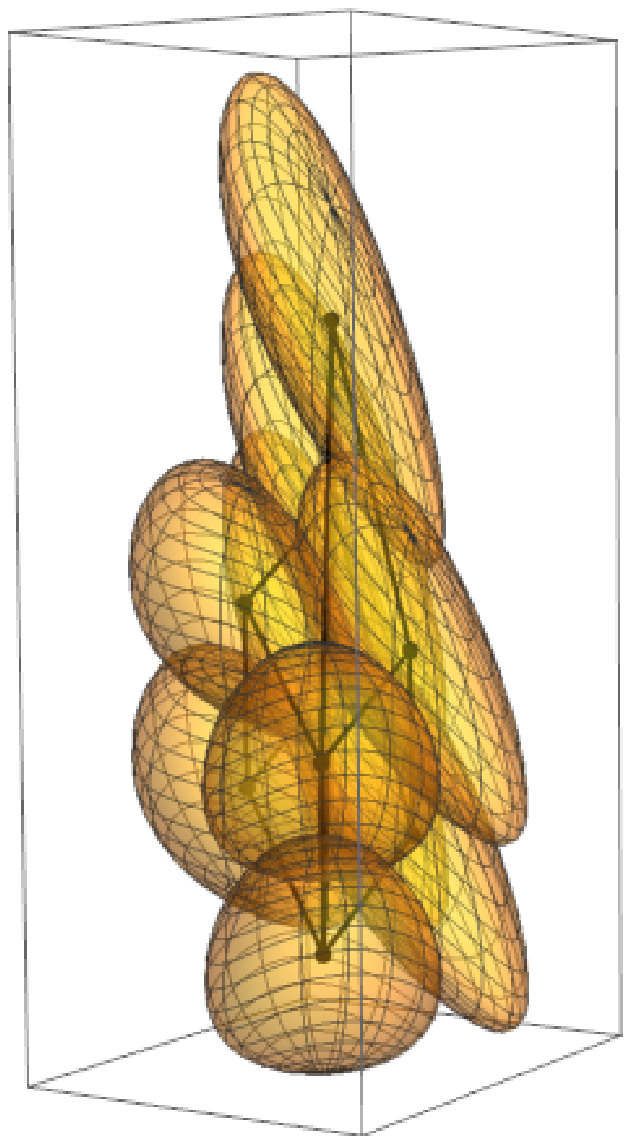}
\caption{Locally optimal lattice-like translation ball covering related to lattice $\Gamma^{u}(\tau_1^u,\tau_2^u,1)$ with density $\approx 1.42783$.} 
\label{}
\end{figure}
The exact determination of the thinnest lattice-like ball covering with congruent translation balls
seems to be difficult, but we are working on refining the upper bound density and 
determine a "good" lower bound density. 

Optimal sphere packings and coverings in other homogeneous Thurston geometries
represent another huge class of open mathematical problems. For
$\NIL$, $\SOL$, $\SLR$, $\HXR$, $\SXR$ geometries only very few results are known
\cite{SchSz11}, \cite{Sz07}, \cite{Sz12-1}, \cite{Sz13-1}, \cite{Sz13-2}, \cite{Sz13-3}. Detailed studies are the objective of
ongoing research.
%

\end{document}